\documentclass[reqno]{amsart}
\usepackage{amsfonts}
\usepackage{amsmath}
\usepackage{amssymb}
\usepackage{graphicx}

\theoremstyle{plain}

\begin{document}
\title[On the Construction of Generalized Bobillier Formula]{On the Construction of Generalized Bobillier Formula}
\author{T\"{u}lay Er\.{i}\d{s}\.{i}r}
\address{T\"{u}lay Er\.{i}\d{s}\.{i}r: Erzincan University, Faculty of Arts and Sciences, Department of Mathematics, Erzincan, Turkey.}
\email{tsoyfidan@sakarya.edu.tr}

\author{Mehmet Al\.{i} G\"{u}ng\"{o}r}
\address{Mehmet Al\.{i} G\"{u}ng\"{o}r: Sakarya University, Faculty of Arts and Sciences, Department of Mathematics, Sakarya,  Turkey.}
\email{agungor@sakarya.edu.tr}

\author{Soley Ersoy}
\address{Soley Ersoy: Sakarya University, Faculty of Arts and Sciences, Department of Mathematics, Sakarya, Turkey.}
\email{sersoy@sakarya.edu.tr}

\date{Received: Month Day, YEAR; Accepted: }
\subjclass[2010]{53A17, 53B50, 11E88.}
\keywords{Bobillier Formula, Euler-Savary Formula, Generalized Complex Numbers.}
\thanks{This paper is in final form and no version of it will be submitted
for publication elsewhere.}

\begin{abstract}
In this study, we consider the generalized complex number system $C_{\rm{p}} = \left\{ {x + iy\;:\;x,y \in R,\;{i^2} = {\rm{p}} \in R} \right\}$ corresponding to elliptical complex number, parabolic complex number and hyperbolic complex number systems for the special cases of $\emph{\emph{p}}<0,  {\kern 1pt} {\kern 1pt} \emph{\emph{p}}=0, {\kern 1pt} {\kern 1pt}  \emph{\emph{p}}<0$, respectively. This system is used to derive Bobillier Formula in the generalized complex plane. In accordance with this purpose we obtain this formula by two different methods for one-parameter planar motion in $C_\emph{\emph{p}}$; the first method depends on using the geometrical interpretation of the generalized Euler-Savary formula and the second one uses the usual relations of the velocities and accelerations.

\end{abstract}

\maketitle

\section{Introduction}

A system of rigid elements (linkages) connected to transmit motion is called mechanism. In other words, a mechanism is a combination of links which can transform a determined motion. A planar mechanism is a mechanical system. So that the trajectories of points in all the bodies of the system is constrained for lying on planes parallel to a grand plane. The rotational axes of hinged joints that connect the bodies in the system are perpendicular to this ground plane, \cite{hunt:kingeo}. Kinematics concerned with the characteristics of movement without considering the concepts of mass and force is sub-branch of mechanics and analyzes the displacement of one point or a point system (object) with respect to time.\\

The Euler-Savary formula which is one of the most used formulas in kinematics and fills an important place in many fields such as mathematics, engineering and astronomy was found by Euler in 1765 and Savary in 1845, \cite{koet:eukin}. This formula giving the relationship between the curvatures of trajectory curves drawn by the points of the moving plane in the fixed plane was studied by \cite{blamu:ebkin,botrot:thekin,buwh:eusav,mull:kinem} in two and three dimensional Euclidean spaces. The Euler-Savary formula in complex plane was given by Masal et al., \cite{matopir:comeu}. Then the Euler-Savary formulas in the hyperbolic plane and Galilean plane were obtained in  \cite{akyu:galeu,erak:hypeul}, respectively. Moreover, Akb{\i}y{\i}k and Y\"{u}ce considered a base curve, a rolling curve and a roulette on complex plane and obtained Euler-Savary formula which gives the relation between the curvatures of these three curves, \cite{yuce}. In addition to these, there are too many studies on the Euler-Savary formula including its applications to kinematics, mechanical engineering of robotics and sciences of machines and mechanisms.\\

The Euler-Savary formulas in Affine Cayley-Klein plane and in generalized complex number plane $ {C_J} = \left\{ {x + Jy{\kern 1pt} {\kern 1pt} {\kern 1pt} {\kern 1pt} :{\kern 1pt} {\kern 1pt} {\kern 1pt} {\kern 1pt} x,y \in R,{\kern 1pt} {\kern 1pt} {\kern 1pt} {\kern 1pt} {\kern 1pt} {\kern 1pt} {\kern 1pt} {J^2} = \emph{\emph{p}}}, {\kern 1pt} {\kern 1pt}{\kern 1pt} {\kern 1pt}{\kern 1pt} {\kern 1pt} \emph{\emph{p}} \in \{-1,0,1\} \right\}\subset {C_\emph{\emph{p}}}$ were given by \cite{nurten} and \cite{nurten2}, respectively.\\

In 1988, Fayet defined a formula giving the relation of the curvatures of second order of one-parameter planar motion and generalizing the Euler-Savary formula. Since this formula analytically solves the problem that the Bobillier's construction solved graphically it is called the Bobillier formula, \cite{fayet:forbob}.\\

In \cite{faye:bobfor}, it was demonstrated that the Bobillier formula can be obtained without the use of the Euler-Savary formula by Fayet. In addition the Bobillier formula was obtained by conventional procedures in \cite{mumi:bovafor}.\\

Ersoy and Bayrak studied the Bobillier formula for the one-parameter planar motion in the complex plane, \cite{erbay:combob}. Also, they investigated the Bobillier formula in Lorentzian sense and saw that the same results can be achieved without the use of the Euler-Savary formula for the Bobillier formula. In doing so, they considered that the cases of polar curves to be timelike or spacelike, separately, \cite{erbay:lobob}.\\

In Galilean plane, G\"{u}rses et al. gave the Bobillier formula by using the geometrical interpretation of Euler-Savary formula in Galilean plane. Moreover, they obtained Bobillier formula without considering the Galilean Euler-Savary formula, \cite{guryuc:galbobi}. Moreover, D\"{u}ndar et al. gave Bobillier formula for the elliptical harmonic motion, \cite{ersoy}.\\

The generalized complex numbers play the same role for Cayley-Klein geometry like that the ordinary numbers play in the Euclidean geometry, \cite{yagl:comnum, yaglom:noneuc}. The Cayley-Klein plane geometries including Euclidean, Galilean, Minkowskian and Bolyai-Lobachevskian were introduced first by F. Klein and A. Cayley, \cite{klein, klein2}. After Cayley and Klein, I. M. Yaglom distinguished these geometries with choosing one of three ways of measuring length (parabolic, elliptic or hyperbolic) between two points on a line and one of the three ways of measuring angles (parabolic, elliptic or hyperbolic) between two lines. This gives nine ways of measuring lengths and angles, \cite{yaglom:noneuc}.\\

In the light of these studies, we construct the Bobillier formula with a new generalization for generalized complex plane $C_\emph{\emph{p}}$. We introduce this generalization in two different ways; using the Euler-Savary formula in the generalized complex plane and an alternative way towards to it.

\section{Preliminaries}
The generalized complex numbers or binary numbers are introduced as follows
$$z = x + iy{\kern 1pt} {\kern 1pt} {\kern 1pt} {\kern 1pt} {\kern 1pt} {\kern 1pt} {\kern 1pt} {\kern 1pt} {\kern 1pt} {\kern 1pt} {\kern 1pt} {\kern 1pt} {\kern 1pt} {\kern 1pt} (x,y \in R),{\kern 1pt} {\kern 1pt} {\kern 1pt} {\kern 1pt} {\kern 1pt} {\kern 1pt} {\kern 1pt} {\kern 1pt} {\kern 1pt} {\kern 1pt} {\kern 1pt} {\kern 1pt} {\kern 1pt} {\kern 1pt} {\kern 1pt} {\kern 1pt} {\kern 1pt} {\kern 1pt} {i^2} = i\emph{\emph{q}} + \emph{\emph{p}}{\kern 1pt} {\kern 1pt} {\kern 1pt} {\kern 1pt} {\kern 1pt} {\kern 1pt} {\kern 1pt} {\kern 1pt} {\kern 1pt} {\kern 1pt} {\kern 1pt} {\kern 1pt} {\kern 1pt} {\kern 1pt} (\emph{\emph{q}},\emph{\emph{p}} \in R).$$
The double, dual and ordinary numbers are the particular members of two parameter family of complex number systems. Moreover, the generalized complex number systems are isomorphic to the double, ordinary and dual complex numbers when $\emph{\emph{p}} + {{{\emph{\emph{q}}^2}} \mathord{\left/
 {\vphantom {{{\emph{\emph{q}}^2}} 4}} \right.
 \kern-\nulldelimiterspace} 4}$  is positive, negative and zero, respectively, \cite{yagl:comnum}.

Unless otherwise stated we assume that ${i^2} = \emph{\emph{p}}$  and  $\emph{\emph{q}} = 0$  $\left( {\emph{\emph{p}} \in R} \right)$. This complex number system is denoted by
$ {C_\emph{\emph{p}}} = \left\{ {x + iy{\kern 1pt} {\kern 1pt} :{\kern 1pt} {\kern 1pt} x,y \in R,{\kern 1pt} {\kern 1pt} {\kern 1pt} {\kern 1pt} {i^2} = \emph{\emph{p}}}, {\kern 1pt} {\kern 1pt}{\kern 1pt} {\kern 1pt} \emph{\emph{p}} \in R \right\}$.

The set ${C_\emph{\emph{p}}}$  is called the generalized complex plane.\\

For ${z_1} = ({x_1} + i{y_1}),{z_2} = ({x_2} + i{y_2}) \in {C_\emph{\emph{p}}}$  the addition, subtraction and product are defined by
$${z_1} \pm {z_2} = ({x_1} + i{y_1}) \pm ({x_2} + i{y_2}) = {x_1} \pm {x_2} + i({y_1} \pm {y_2})$$
and
$${M^\emph{\emph{p}}}({z_1},{z_2}) = ({x_1}{x_2} + \emph{\emph{p}}{y_1}{y_2}) + i({x_1}{y_2} + {x_2}{y_1}),$$
respectively. The product definition yields the ordinary, Study and Clifford products as  $\emph{\emph{p}}$ is equal to  $ - 1$, $0$ and  $1$
$${M^{ - 1}}({z_1},{z_2}) = ({x_1}{x_2} - {y_1}{y_2}) + i({x_1}{y_2} + {x_2}{y_1}){\kern 1pt} {\kern 1pt} {\kern 1pt} {\kern 1pt}   {\kern 1pt} {\kern 1pt} \emph{\emph{for}}{\kern 1pt} {\kern 1pt} {\kern 1pt} {\kern 1pt} {\kern 1pt} {\kern 1pt}  \emph{\emph{p}} =  - 1,$$
$${M^0}({z_1},{z_2}) = ({x_1}{x_2}) + i({x_1}{y_2} + {x_2}{y_1}) {\kern 1pt} {\kern 1pt} {\kern 1pt} {\kern 1pt}  {\kern 1pt} {\kern 1pt}     \emph{\emph{for}}{\kern 1pt} {\kern 1pt} {\kern 1pt} {\kern 1pt} {\kern 1pt} {\kern 1pt}  \emph{\emph{p}} = 0,$$
$${M^1}({z_1},{z_2}) = ({x_1}{x_2} + {y_1}{y_2}) + i({x_1}{y_2} + {x_2}{y_1}) {\kern 1pt}{\kern 1pt} {\kern 1pt} {\kern 1pt} {\kern 1pt}  {\kern 1pt}    \emph{\emph{for}} {\kern 1pt} {\kern 1pt} {\kern 1pt} {\kern 1pt} {\kern 1pt} {\kern 1pt}  \emph{\emph{p}} = 1,$$
respectively, \cite{yagl:comnum}.

For  ${\textbf{z}_1} = ({x_1},{y_1}),{\kern 1pt} {\kern 1pt} {\kern 1pt} {\kern 1pt} {\kern 1pt} {\kern 1pt} {\kern 1pt} {\kern 1pt} {\textbf{z}_2} = ({x_2},{y_2}) \in {C_\emph{\emph{p}}}$, the scalar product is as follows
\begin{equation}
{\left\langle {{\textbf{z}_1},{\textbf{z}_2}} \right\rangle _\emph{\emph{p}}} = {\mathop{\rm Re}\nolimits} \left( {{M^\emph{\emph{p}}}({\textbf{z}_1},\overline {{\textbf{z}_2}} )} \right) = {\mathop{\rm Re}\nolimits} \left( {{M^\emph{\emph{p}}}(\overline {{\textbf{z}_1}} ,{\textbf{z}_2})} \right) = {x_1}{x_2} - \emph{\emph{p}}{y_1}{y_2}
\end{equation}
and the  magnitude of $\textbf{z} = (x,y) \in {C_\emph{\emph{p}}}$  is given by the non-negative real number
\begin{equation}
{\left\| \textbf{z} \right\|_\emph{\emph{p}}} = \sqrt {\left| {{M^\emph{\emph{p}}}(\textbf{z}, \bar{\textbf{z}} )} \right|}  = \sqrt {\left| {{x^2} - \emph{\emph{p}}{y^2}} \right|}
\end{equation}
where the over bar denotes the complex conjugation. The equalities
\begin{equation}
{\left\langle {{\textbf{z}_1},{\textbf{z}_2}} \right\rangle _\emph{\emph{p}}} = {\left\| {{\textbf{z}_1}} \right\|_\emph{\emph{p}}}{\left\| {{\textbf{z}_2}} \right\|_\emph{\emph{p}}}\cos \emph{\emph{p}}{\theta _\emph{\emph{p}}},
\end{equation}
and for $\emph{\emph{p}} \ne 0$
\begin{equation}
{\left\| {{\textbf{z}_1}{ \wedge _\emph{\emph{p}}}{\textbf{z}_2}} \right\|_\emph{\emph{p}}} = \sqrt {\left| { - \emph{\emph{p}}} \right|} {\left\| {{\textbf{z}_1}} \right\|_\emph{\emph{p}}}{\left\| {{\textbf{z}_2}} \right\|_\emph{\emph{p}}}\sin \emph{\emph{p}}{\theta _\emph{\emph{p}}}
\end{equation}
and for $\emph{\emph{p}} = 0$, (from the special definition in Galilean plane)
\begin{equation}
{\left\| {{\textbf{z}_1}{ \wedge _\emph{\emph{p}}}{\textbf{z}_2}} \right\|_\emph{\emph{p}}} ={\left\| {{\textbf{z}_1}} \right\|_\emph{\emph{p}}}{\left\| {{\textbf{z}_2}} \right\|_\emph{\emph{p}}}\sin \emph{\emph{p}}{\theta _\emph{\emph{p}}}
\end{equation}
are valid in the generalized complex plane where ${\theta _\emph{\emph{p}}}$  is the  $\emph{\emph{p}}-$rotation angle between the vectors ${\textbf{z}_1}$  and  ${\textbf{z}_2}$.\\

In the Euclidean plane, a circle is defined in two ways. The first definition is the locus of points a fixed distance from a given fixed point. The second is the set of points from which a segment $AB$  is seen at a constant directed angle  ${\theta _\emph{\emph{p}}}$. Using these definitions, since a unit circle is the geometric locus of points  $z$ satisfying  ${\left\| z \right\|_\emph{\emph{p}}} = 1$, the unit circles in the  generalized complex planes are drawn as in Figure 1.

\begin{center}
\includegraphics[scale=0.75]{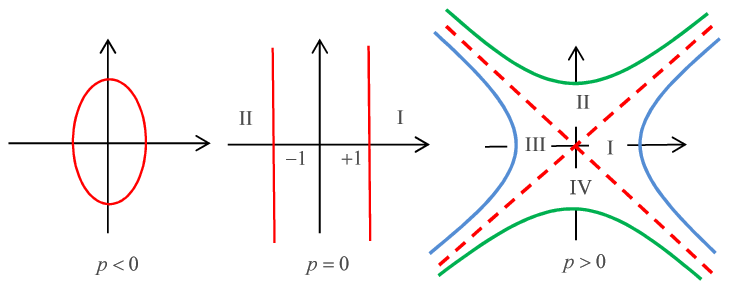}\

\scriptsize
{Figure 1. Unit circles in $C_\emph{\emph{p}}$.}
\normalsize
\end{center}

Let us interpret Figure 1. If we take  $\emph{\emph{p}} < 0$, then the unit circle corresponds to unit ellipse of the form ${x^2} + \left| \emph{\emph{p}} \right|{y^2} = 1$  and the complex number system  ${C_\emph{\emph{p}}}\left( {\emph{\emph{p}} < 0} \right)$ corresponds to elliptical complex number system. In this case, if  $\emph{\emph{p}} =  - 1$, the unit circle in the  generalized complex plane corresponds to a standard unit circle defined by ${x^2} + {y^2} = 1$  and the  generalized complex plane corresponds to the Euclidean plane. \\

If we consider  $\emph{\emph{p}} = 0$, the above \textbf{\emph{\emph{circle}}} definitions yield different sets of points in Galilean geometry. According to the first definition, the equation  $\left\| z \right\|_0^2 = {x^2}$ is hold and the unit \textbf{\emph{\emph{circle}}} is as  $x =  \pm 1$, (Figure 1) and for the first circle we cannot talk of curvature. According to the second definition, the set of points is called a \textbf{\emph{\emph{cycle}}}.  The \textbf{\emph{\emph{cycle}}} (inflection cycle, red lines) on Galilean plane is the (Euclidean) parabola, \cite{yaglom:noneuc}. Therefore, the space ${C_0}$  is the parabolic complex number system and the  generalized complex plane corresponds to the Galilean plane. Moreover, as can be seen from Figure 1 the parabolic complex plane is divided into two area by the imaginary axis. \\

Finally, when  $\emph{\emph{p}} > 0$, unit circles in ${C_\emph{\emph{p}}}\left( {\emph{\emph{p}} > 0} \right)$  are the hyperbolas of the form $\left| {{x^2} - \emph{\emph{p}}{y^2}} \right| = 1$  whose asymptotes are  $y =  \pm {x \mathord{\left/
 {\vphantom {x {\sqrt \emph{\emph{p}} }}} \right.
 \kern-\nulldelimiterspace} {\sqrt \emph{\emph{p}} }}$ (red dashed line in Figure 1). Thus, the space  ${C_\emph{\emph{p}}}\left( {\emph{\emph{p}} > 0} \right)$ corresponds to hyperbolic complex number system. Especially, if we take  $\emph{\emph{p}} = 1$, the  generalized complex plane is referred to as a well-known hyperbolic plane and the asymptotes of the unit circles separate the hyperbolic planes into four regions, \cite{hahar:genkom}.\\

$\emph{\emph{p}}-$trigonometric functions ( $\emph{\emph{p}}-$cosine ($\cos \emph{\emph{p}}$),  $\emph{\emph{p}}-$sine ($\sin \emph{\emph{p}}$) and  $\emph{\emph{p}}-$tangent ($\tan \emph{\emph{p}}$) are, respectively, defined by
\begin{equation}
\begin{array}{l}
\cos \emph{\emph{p}}{\theta _\emph{\emph{p}}} = \left\{ \begin{array}{l}
\cos \left( {{\theta _\emph{\emph{p}}}\sqrt {\left| \emph{\emph{p}} \right|} } \right),{\kern 1pt} {\kern 1pt} {\kern 1pt} {\kern 1pt} {\kern 1pt} {\kern 1pt} {\kern 1pt} {\kern 1pt} {\kern 1pt} {\kern 1pt} {\kern 1pt} {\kern 1pt} {\kern 1pt} {\kern 1pt} {\kern 1pt} {\kern 1pt} {\kern 1pt} {\kern 1pt} {\kern 1pt} {\kern 1pt} {\kern 1pt} {\kern 1pt} {\kern 1pt} \emph{\emph{p}} < 0\\
1,{\kern 1pt} {\kern 1pt} {\kern 1pt} {\kern 1pt} {\kern 1pt} {\kern 1pt} {\kern 1pt} {\kern 1pt} {\kern 1pt} {\kern 1pt} {\kern 1pt} {\kern 1pt} {\kern 1pt} {\kern 1pt} {\kern 1pt} {\kern 1pt} {\kern 1pt} {\kern 1pt} {\kern 1pt} {\kern 1pt} {\kern 1pt} {\kern 1pt} {\kern 1pt} {\kern 1pt} {\kern 1pt} {\kern 1pt} {\kern 1pt} {\kern 1pt} {\kern 1pt} {\kern 1pt} {\kern 1pt} {\kern 1pt} {\kern 1pt} {\kern 1pt} {\kern 1pt} {\kern 1pt} {\kern 1pt} {\kern 1pt} {\kern 1pt} {\kern 1pt} {\kern 1pt} {\kern 1pt} {\kern 1pt} {\kern 1pt} {\kern 1pt} {\kern 1pt} {\kern 1pt} {\kern 1pt} {\kern 1pt} {\kern 1pt} {\kern 1pt} {\kern 1pt} {\kern 1pt} {\kern 1pt} {\kern 1pt} {\kern 1pt} {\kern 1pt} {\kern 1pt} {\kern 1pt} {\kern 1pt} {\kern 1pt} {\kern 1pt} {\kern 1pt} {\kern 1pt} {\kern 1pt} {\kern 1pt} {\kern 1pt} {\kern 1pt} {\kern 1pt} {\kern 1pt} {\kern 1pt} {\kern 1pt} {\kern 1pt} {\kern 1pt} {\kern 1pt} {\kern 1pt} {\kern 1pt} {\kern 1pt} {\kern 1pt} \emph{\emph{p}} = 0{\kern 1pt} {\kern 1pt} {\kern 1pt} {\kern 1pt} {\kern 1pt} {\kern 1pt} {\kern 1pt} (\emph{branch}{\kern 1pt} {\kern 1pt} I)\\
\cosh \left( {{\theta _\emph{\emph{p}}}\sqrt \emph{\emph{p}} } \right),{\kern 1pt} {\kern 1pt} {\kern 1pt} {\kern 1pt} {\kern 1pt} {\kern 1pt} {\kern 1pt} {\kern 1pt} {\kern 1pt} {\kern 1pt} {\kern 1pt} {\kern 1pt} {\kern 1pt} {\kern 1pt} {\kern 1pt} {\kern 1pt} {\kern 1pt} {\kern 1pt} {\kern 1pt} {\kern 1pt} {\kern 1pt} {\kern 1pt}{\kern 1pt} {\kern 1pt}{\kern 1pt} {\kern 1pt}{\kern 1pt} {\kern 1pt}\emph{\emph{p}} > 0{\kern 1pt} {\kern 1pt} {\kern 1pt} {\kern 1pt} {\kern 1pt} {\kern 1pt} {\kern 1pt} {\kern 1pt} {\kern 1pt} (\emph{branch}{\kern 1pt} {\kern 1pt} II)
\end{array} \right.,{\kern 1pt} {\kern 1pt} {\kern 1pt} {\kern 1pt} {\kern 1pt} {\kern 1pt} {\kern 1pt} {\kern 1pt} \\
\sin \emph{\emph{p}}{\theta _\emph{\emph{p}}} = \left\{ \begin{array}{l}
\frac{1}{{\sqrt {\left| \emph{\emph{p}} \right|} }}\sin \left( {{\theta _\emph{\emph{p}}}\sqrt {\left| \emph{\emph{p}} \right|} } \right),{\kern 1pt} {\kern 1pt} {\kern 1pt} {\kern 1pt} {\kern 1pt} {\kern 1pt} \emph{\emph{p}} < 0\\
{\theta _\emph{\emph{p}}},{\kern 1pt} {\kern 1pt} {\kern 1pt} {\kern 1pt} {\kern 1pt} {\kern 1pt} {\kern 1pt} {\kern 1pt} {\kern 1pt} {\kern 1pt} {\kern 1pt} {\kern 1pt} {\kern 1pt} {\kern 1pt} {\kern 1pt} {\kern 1pt} {\kern 1pt} {\kern 1pt} {\kern 1pt} {\kern 1pt} {\kern 1pt} {\kern 1pt} {\kern 1pt} {\kern 1pt} {\kern 1pt} {\kern 1pt} {\kern 1pt} {\kern 1pt} {\kern 1pt} {\kern 1pt} {\kern 1pt} {\kern 1pt} {\kern 1pt} {\kern 1pt} {\kern 1pt} {\kern 1pt} {\kern 1pt} {\kern 1pt} {\kern 1pt} {\kern 1pt} {\kern 1pt} {\kern 1pt} {\kern 1pt} {\kern 1pt} {\kern 1pt} {\kern 1pt} {\kern 1pt} {\kern 1pt} {\kern 1pt} {\kern 1pt} {\kern 1pt} {\kern 1pt} {\kern 1pt} {\kern 1pt} {\kern 1pt} {\kern 1pt} {\kern 1pt} {\kern 1pt} {\kern 1pt} {\kern 1pt} {\kern 1pt} {\kern 1pt} {\kern 1pt} {\kern 1pt} {\kern 1pt} {\kern 1pt} {\kern 1pt} {\kern 1pt} {\kern 1pt} {\kern 1pt} {\kern 1pt} {\kern 1pt} {\kern 1pt} {\kern 1pt} {\kern 1pt} {\kern 1pt} {\kern 1pt} \emph{\emph{p}} = 0{\kern 1pt} {\kern 1pt} {\kern 1pt} {\kern 1pt} {\kern 1pt} {\kern 1pt} {\kern 1pt} {\kern 1pt} {\kern 1pt} (\emph{branch}{\kern 1pt} {\kern 1pt} I)\\
\frac{1}{{\sqrt \emph{\emph{p}} }}\sinh \left( {{\theta _\emph{\emph{p}}}\sqrt \emph{\emph{p}} } \right),{\kern 1pt} {\kern 1pt} {\kern 1pt} {\kern 1pt} {\kern 1pt} {\kern 1pt} {\kern 1pt} {\kern 1pt} {\kern 1pt} {\kern 1pt}{\kern 1pt} {\kern 1pt}{\kern 1pt} {\kern 1pt}{\kern 1pt} {\kern 1pt}{\kern 1pt} {\kern 1pt}\emph{\emph{p}} > 0\,\,{\kern 1pt} {\kern 1pt} {\kern 1pt} {\kern 1pt} {\kern 1pt} {\kern 1pt} (\emph{branch}{\kern 1pt} {\kern 1pt} II)
\end{array} \right.,{\kern 1pt} {\kern 1pt} {\kern 1pt} {\kern 1pt} {\kern 1pt} {\kern 1pt} {\kern 1pt} {\kern 1pt} {\kern 1pt} {\kern 1pt} {\kern 1pt} {\kern 1pt} {\kern 1pt} {\kern 1pt} {\kern 1pt} {\kern 1pt} {\kern 1pt} {\kern 1pt}
\end{array}
\end{equation}
and
\begin{equation}
\tan \emph{\emph{p}}{\theta _\emph{\emph{p}}} = \frac{{\sin \emph{\emph{p}}{\theta _\emph{\emph{p}}}}}{{\cos \emph{\emph{p}}{\theta _\emph{\emph{p}}}}}.
\end{equation}
Additionally, the derivatives of  $\emph{\emph{p}}-$trigonometric functions are as follows,
\begin{equation}
\frac{d}{{d{\theta _\emph{\emph{p}}}}}(\cos \emph{\emph{p}}{\theta _\emph{\emph{p}}}) = \emph{\emph{p}}\sin \emph{\emph{p}}{\theta _\emph{\emph{p}}},\,\,\,\,\,\,\,\,\,\,\,\,\,\,\,\,\,\,\,\,\,\frac{d}{{d{\theta _\emph{\emph{p}}}}}(\sin \emph{\emph{p}}{\theta _\emph{\emph{p}}}) = \cos \emph{\emph{p}}{\theta _\emph{\emph{p}}}
\end{equation}
and the generalized Euler's formula is as follows
\begin{equation}
{e^{i{\theta _\emph{\emph{p}}}}} = \cos \emph{\emph{p}}{\theta _\emph{\emph{p}}} + i\sin \emph{\emph{p}}{\theta _\emph{\emph{p}}}.
\end{equation}
The polar and exponential forms of any  generalized complex number $z$  are
\[z = {r_\emph{\emph{p}}}(\cos\emph{\emph{p}}{\theta _\emph{\emph{p}}} + i\sin \emph{\emph{p}}{\theta _\emph{\emph{p}}}) = {r_\emph{\emph{p}}}{e^{i{\theta _\emph{\emph{p}}}}}\]
where ${r_\emph{\emph{p}}} = {\left\| z \right\|_\emph{\emph{p}}}$  and ${\theta _\emph{\emph{p}}}$  are the  $\emph{\emph{p}}-$magnitude and  argument of  $z$, respectively, \cite{hahar:genkom}.\\

Let $K_\emph{\emph{p}}$  and $K'_\emph{\emph{p}}$  be the moving and fixed  generalized complex planes and $\left\{ {O;{\textbf{t}_1},\,{\textbf{t}_2}} \right\}$  and $\left\{ {O';{{\textbf{t}'}_1},\,{{\textbf{t}'}_2}} \right\}$  be the perpendicular coordinate system of these planes, respectively. Moreover, we consider that the one-parameter planar motion $K_\emph{\emph{p}}/K'_\emph{\emph{p}}$ in the generalized complex plane $C_\emph{\emph{p}}$. The geometric loci of the pole points in the fixed and moving  generalized complex planes $K'_\emph{\emph{p}}$ and $K_\emph{\emph{p}}$  are called fixed pole curve $(Q')$  and moving pole curve  $(Q)$, respectively. \\

A point $X$ taken in the moving  generalized complex plane $K_\emph{\emph{p}}$ draws a trajectory which has curvature center $X'$  in the fixed  generalized complex plane  $K'_\emph{\emph{p}}$. Conversely, in the reverse motion $K'_\emph{\emph{p}}/K_\emph{\emph{p}}$ in $C_\emph{\emph{p}}$ the point $X'$  in  $K'_\emph{\emph{p}}$ draws a trajectory with curvature center $X$  in  $K_\emph{\emph{p}}$. This relation between the points $X$  and  $X'$ is given by Euler-Savary Formula with the following equation

\begin{equation}\label{10}
\frac{1}{{r'}} - \frac{1}{r} = {\mathop{\rm Im}\nolimits} \left( {{e^{i\theta _{_\emph{\emph{p}}}}}} \right)\left( {\frac{1}{{{a'}}} - \frac{1}{{{a}}}} \right)
\end{equation}
where $r$  and $r'$  are the radii of curvature of the pole curves  $(Q)$ and  $(Q')$ of one-parameter planar motion $K_\emph{\emph{p}}/K'_\emph{\emph{p}}$ and $a$, $a'$ represent the distances from the rotation pole to the point $X,X'$, respectively. This formula for homothetic motion on the plane $ {C_J} = \left\{ {x + Jy {\kern 1pt} :{\kern 1pt} x,y \in R,{\kern 1pt} {\kern 1pt} {J^2} = \emph{\emph{p}}},{\kern 1pt} {\kern 1pt} \emph{\emph{p}} \in \{-1,0,1\} \right\}\subset {C_\emph{\emph{p}}}$ was given by G\"{u}rses et al.

\section{On the Construction of Generalized Bobillier Formula}

In this section, we use a method to derive generalized Bobillier formula for generalized complex plane $C_\emph{\emph{p}}$. For this, we regard to the generalized Euler-Savary formula given by the equation (\ref{10}).

Let $K_\emph{\emph{p}}$  and $K'_\emph{\emph{p}}$  be the moving and fixed  generalized complex planes and $\left\{ {O;{\textbf{t}_1},\,{\textbf{t}_2}} \right\}$  and $\left\{ {O';{{\textbf{t}'}_1},\,{{\textbf{t}'}_2}} \right\}$  be the perpendicular coordinate system of these planes, respectively. Moreover, we assume that the points $N_\emph{\emph{p}}^1,$ $N_\emph{\emph{p}}^2$  and $N_\emph{\emph{p}}^3$ are the fixed points in the moving generalized complex plane $K_\emph{\emph{p}}$  and the points $\gamma _\emph{\emph{p}}^1,$  $\gamma _\emph{\emph{p}}^2$ and $\gamma _\emph{\emph{p}}^3$  are centers of curvature of trajectory drawn by these fixed points in the fixed  generalized complex plane  $K'_\emph{\emph{p}}$.

The normals of trajectories drawn by these points pass from the instantaneous rotation center $I$  and there is a  $\emph{\emph{p}}-$rotation pole at each  $t$ moment in $C_\emph{\emph{p}}$. The pole curves rolls upon each other without sliding during the motion of the generalized complex planes  $K_\emph{\emph{p}}/K'_\emph{\emph{p}}$. So, the pole curves  $(Q)$, $(Q')$  are tangent to each other and have the same velocity at each  $t$ moment. Thus, the real axis is the common tangent and the imaginary axis is the common normal for the pole curves in $C_\emph{\emph{p}}$.

If ${\theta _\emph{\emph{p}}}$  is the  $\emph{\emph{p}}-$rotation angle of motion of the plane with respect to $K'_\emph{\emph{p}}$  at each $t$  moment, then each point $N_\emph{\emph{p}}^i$  makes a rotation motion with ${\dot \theta _\emph{\emph{p}}}$  angular velocity at the instantaneous center  $I$.

Let
\begin{equation}
\textbf{X}_\emph{\emph{p}}^1 = \frac{{\textbf{IN}_\emph{\emph{p}}^1}}{{\left\| {\textbf{IN}_\emph{\emph{p}}^1} \right\|}},{\kern 1pt} {\kern 1pt} {\kern 1pt} {\kern 1pt} {\kern 1pt} {\kern 1pt} {\kern 1pt} {\kern 1pt} {\kern 1pt} {\kern 1pt} {\kern 1pt} {\kern 1pt} {\kern 1pt} {\kern 1pt} {\kern 1pt} {\kern 1pt} {\kern 1pt} {\kern 1pt} {\kern 1pt} {\kern 1pt} {\kern 1pt} {\kern 1pt} {\kern 1pt} {\kern 1pt} {\kern 1pt} {\kern 1pt} {\kern 1pt} \textbf{X}_\emph{\emph{p}}^2 = \frac{{\textbf{IN}_\emph{\emph{p}}^2}}{{\left\| {\textbf{IN}_\emph{\emph{p}}^2} \right\|}},{\kern 1pt} {\kern 1pt} {\kern 1pt} {\kern 1pt} {\kern 1pt} {\kern 1pt} {\kern 1pt} {\kern 1pt} {\kern 1pt} {\kern 1pt} {\kern 1pt} {\kern 1pt} {\kern 1pt} {\kern 1pt} {\kern 1pt} {\kern 1pt} {\kern 1pt} {\kern 1pt} {\kern 1pt} {\kern 1pt} {\kern 1pt} {\kern 1pt} {\kern 1pt} {\kern 1pt} {\kern 1pt} {\kern 1pt} \textbf{X}_\emph{\emph{p}}^3 = \frac{{\textbf{IN}_\emph{\emph{p}}^3}}{{\left\| {\textbf{IN}_\emph{\emph{p}}^3} \right\|}}
\end{equation}
be the unit vectors in the direction of the pole rays $\textbf{IN}_\emph{\emph{p}}^1,$  $\textbf{IN}_\emph{\emph{p}}^2$  and  $\textbf{IN}_\emph{\emph{p}}^3$, respectively, (see Figure 2 for the special cases of $\emph{\emph{p}}=-1,0,+1$).

	The  $\emph{\emph{p}}-$distances of the points  $N_\emph{\emph{p}}^k$ and $\gamma _\emph{\emph{p}}^k$  from the origin  $I$ are ${\rho _k}$  and    ${\rho '_k}$ $(k = 1,2,3)$, respectively, then from the equations (3) and (4) we can write
\begin{equation}
\left| {{{\left\langle {\textbf{IN}_\emph{\emph{p}}^1,\textbf{X}_\emph{\emph{p}}^1} \right\rangle }_\emph{\emph{p}}}} \right| = {\rho _1},{\kern 1pt} {\kern 1pt} {\kern 1pt} {\kern 1pt} {\kern 1pt} {\kern 1pt} {\kern 1pt} {\kern 1pt} {\kern 1pt} {\kern 1pt} {\kern 1pt} {\kern 1pt} {\kern 1pt} {\kern 1pt} {\kern 1pt} {\kern 1pt} {\kern 1pt} {\kern 1pt} {\kern 1pt} {\kern 1pt} {\kern 1pt} \left| {{{\left\langle {\textbf{I}\gamma _\emph{\emph{p}}^1,\textbf{X}_\emph{\emph{p}}^1} \right\rangle }_\emph{\emph{p}}}} \right| = {\rho '_1}.
\end{equation}
Similarly,
$$\left| {{{\left\langle {\textbf{IN}_\emph{\emph{p}}^2,\textbf{X}_\emph{\emph{p}}^2} \right\rangle }_\emph{\emph{p}}}} \right| = {\rho _2},{\kern 1pt} {\kern 1pt} {\kern 1pt} {\kern 1pt} {\kern 1pt} {\kern 1pt} {\kern 1pt} {\kern 1pt} {\kern 1pt} {\kern 1pt} {\kern 1pt} {\kern 1pt} {\kern 1pt} {\kern 1pt} {\kern 1pt} {\kern 1pt} {\kern 1pt} {\kern 1pt} {\kern 1pt} {\kern 1pt} {\kern 1pt} {\kern 1pt} {\kern 1pt} {\kern 1pt} \left| {{{\left\langle {\textbf{I}\gamma _\emph{\emph{p}}^2,\textbf{X}_\emph{\emph{p}}^2} \right\rangle }_\emph{\emph{p}}}} \right| = {\rho '_2}$$
and
$$\left| {{{\left\langle {\textbf{IN}_\emph{\emph{p}}^3,\textbf{X}_\emph{\emph{p}}^3} \right\rangle }_\emph{\emph{p}}}} \right| = {\rho _3},{\kern 1pt} {\kern 1pt} {\kern 1pt} {\kern 1pt} {\kern 1pt} {\kern 1pt} {\kern 1pt} {\kern 1pt} {\kern 1pt} {\kern 1pt} {\kern 1pt} {\kern 1pt} {\kern 1pt} {\kern 1pt} {\kern 1pt} {\kern 1pt} {\kern 1pt} {\kern 1pt} {\kern 1pt} {\kern 1pt} {\kern 1pt} {\kern 1pt} {\kern 1pt} {\kern 1pt} \left| {{{\left\langle {\textbf{I}\gamma _\emph{\emph{p}}^3,\textbf{X}_\emph{\emph{p}}^3} \right\rangle }_\emph{\emph{p}}}} \right| = {\rho '_3}.$$

	An inflection point may be defined to be a point whose trajectory momentarily has an infinite radius of curvature. Such points also have zero acceleration normal to their trajectory, \cite{blamu:ebkin, saerhurag:eusaveq, saarra:doubval}.  Let the inflection points be $N{_\emph{\emph{p}}^1}^*,$  $N{_\emph{\emph{p}}^2}^*$  and   $N{_\emph{\emph{p}}^3}^*$. The locus of such points is a circle in the moving  generalized complex plane $K_\emph{\emph{p}}$  called as an inflection circle. So, from the equation (3) and (4) we can write
\scriptsize
\begin{equation}
\left| {{{\left\langle {\textbf{IN}{{_\emph{\emph{p}}^1}^*},\textbf{X}_\emph{\emph{p}}^1} \right\rangle }_\emph{\emph{p}}}} \right| = \rho _1^*,{\kern 1pt} {\kern 1pt} {\kern 1pt} {\kern 1pt} {\kern 1pt} {\kern 1pt} {\kern 1pt} {\kern 1pt} {\kern 1pt} {\kern 1pt} {\kern 1pt} {\kern 1pt} {\kern 1pt} {\kern 1pt} {\kern 1pt} {\kern 1pt} {\kern 1pt} {\kern 1pt} {\kern 1pt} {\kern 1pt} \left| {{{\left\langle {\textbf{IN}{{_\emph{\emph{p}}^2}^*},\textbf{X}_\emph{\emph{p}}^2} \right\rangle }_\emph{\emph{p}}}} \right|{\kern 1pt}  = \rho _2^*,{\kern 1pt} {\kern 1pt} {\kern 1pt} {\kern 1pt} {\kern 1pt} {\kern 1pt} {\kern 1pt} {\kern 1pt} {\kern 1pt} {\kern 1pt} {\kern 1pt} {\kern 1pt} {\kern 1pt} {\kern 1pt} {\kern 1pt} {\kern 1pt} {\kern 1pt} {\kern 1pt} {\kern 1pt} {\kern 1pt} \left| {{{\left\langle {\textbf{IN}{{_\emph{\emph{p}}^3}^*},\textbf{X}_\emph{\emph{p}}^3} \right\rangle }_\emph{\emph{p}}}} \right|{\kern 1pt}  = \rho _3^*
\end{equation}

\begin{center}
\includegraphics[scale=0.8]{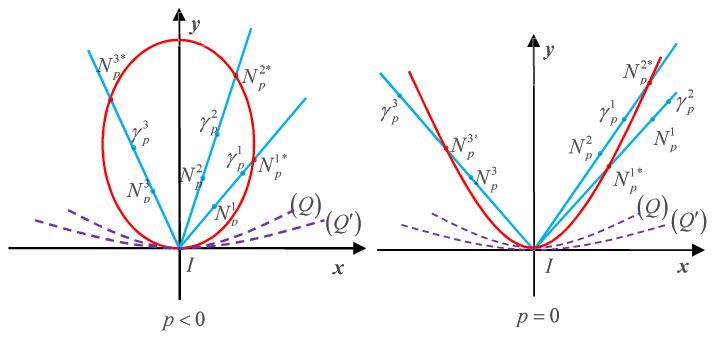}\

\includegraphics[scale=0.8]{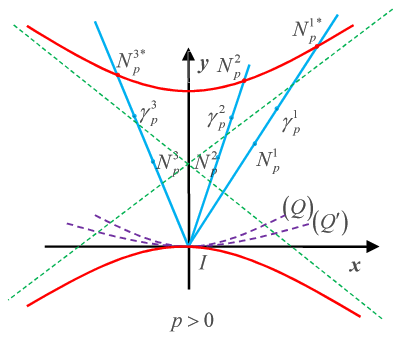}\

\scriptsize
{Figure 2. The special cases of the inflection circles of motion in the generalized complex plane.}
\end{center}
\normalsize

On the other hand the points move around an trajectory whose instantaneous center is $I$  during the  motion  $K_\emph{\emph{p}}/K'_\emph{\emph{p}}$. So, the amount of displacement of the point $I$  is equal to the product of the diameter $h$  with the amount of angular displacement.\\

We take that the inflection points  $N{_\emph{\emph{p}}^1}^*,$ $N{_\emph{\emph{p}}^2}^*$   and $N{_\emph{\emph{p}}^3}^*$  are on the direction of $(I,\textbf{X}_\emph{\emph{p}}^1),$  $(I,\textbf{X}_\emph{\emph{p}}^2)$  and  $(I,\textbf{X}_\emph{\emph{p}}^3)$, respectively. Moreover, we assume that, the images of these inflection points are  $Q_\emph{\emph{p}}^1,$  $Q_\emph{\emph{p}}^2$ and $Q_\emph{\emph{p}}^3$  where  $\textbf{IQ}_\emph{\emph{p}}^k = \frac{1}{{\rho _k^*}}\textbf{X}_\emph{\emph{p}}^k{\kern 1pt} {\kern 1pt} {\kern 1pt} ,{\kern 1pt} {\kern 1pt} {\kern 1pt} {\kern 1pt} {\kern 1pt} {\kern 1pt} {\kern 1pt} \left( {1 \le k \le 3} \right)$, (Figure 3,4,5).  Then, from the equations (3) and (4) we have
\scriptsize
\begin{equation}
\left| {{{\left\langle {\textbf{IQ}_\emph{\emph{p}}^1,\textbf{X}_\emph{\emph{p}}^1} \right\rangle }_\emph{\emph{p}}}} \right| = \frac{1}{{\rho _1^*}},{\kern 1pt} {\kern 1pt} {\kern 1pt} {\kern 1pt} {\kern 1pt} {\kern 1pt} {\kern 1pt} {\kern 1pt} {\kern 1pt} {\kern 1pt} {\kern 1pt} {\kern 1pt} {\kern 1pt} \left| {{{\left\langle {\textbf{IQ}_\emph{\emph{p}}^2,\textbf{X}_\emph{\emph{p}}^2} \right\rangle }_\emph{\emph{p}}}} \right| = \frac{1}{{\rho _2^*}},{\kern 1pt} {\kern 1pt} {\kern 1pt} {\kern 1pt} {\kern 1pt} {\kern 1pt} {\kern 1pt} {\kern 1pt} {\kern 1pt} {\kern 1pt} {\kern 1pt} \left| {{{\left\langle {\textbf{IQ}_\emph{\emph{p}}^3,\textbf{X}_\emph{\emph{p}}^3} \right\rangle }_\emph{\emph{p}}}} \right| = \frac{1}{{\rho _3^*}}.
\end{equation}
\normalsize
So, the following equations hold;
\begin{equation}
\begin{array}{l}
\textbf{IQ}_\emph{\emph{p}}^1\sin \emph{\emph{p}}\theta _\emph{\emph{p}}^1 = \frac{1}{{\rho _1^*}}\textbf{X}_\emph{\emph{p}}^1\sin \emph{\emph{p}}\theta _\emph{\emph{p}}^1 = \frac{1}{h}\textbf{X}_\emph{\emph{p}}^1,{\kern 1pt} {\kern 1pt} {\kern 1pt} {\kern 1pt} {\kern 1pt} {\kern 1pt} {\kern 1pt} {\kern 1pt} {\kern 1pt} \\
\textbf{IQ}_\emph{\emph{p}}^2\sin \emph{\emph{p}}\theta _\emph{\emph{p}}^2 = \frac{1}{{\rho _2^*}}\textbf{X}_p\emph{\emph{p}}^2\sin \emph{\emph{p}}\theta _\emph{\emph{p}}^2 = \frac{1}{h}\textbf{X}_\emph{\emph{p}}^2,{\kern 1pt} \\
{\kern 1pt} \textbf{IQ}_\emph{\emph{p}}^3\sin \emph{\emph{p}}\theta _\emph{\emph{p}}^3 = \frac{1}{{\rho _3^*}}\textbf{X}_\emph{\emph{p}}^3\sin \emph{\emph{p}}\theta _\emph{\emph{p}}^3 = \frac{1}{h}\textbf{X}_\emph{\emph{p}}^3.
\end{array}
\end{equation}
where the diameter of the inflection cycle is $h$.
The last three equations show that
\begin{equation}
\left| {{{\left\langle {\textbf{IQ}_\emph{\emph{p}}^1,\textbf{X}_\emph{\emph{p}}^1} \right\rangle }_\emph{\emph{p}}}} \right|\sin \emph{\emph{p}}\theta _\emph{\emph{p}}^1 = \left| {{{\left\langle {\textbf{IQ}_\emph{\emph{p}}^2,\textbf{X}_\emph{\emph{p}}^2} \right\rangle }_\emph{\emph{p}}}} \right|\sin \emph{\emph{p}}\theta _\emph{\emph{p}}^2 = \left| {{{\left\langle {\textbf{IQ}_\emph{\emph{p}}^3,\textbf{X}_\emph{\emph{p}}^3} \right\rangle }_\emph{\emph{p}}}} \right|\sin \emph{\emph{p}}\theta _\emph{\emph{p}}^3 = \frac{1}{h}.
\end{equation}

So, the set of the points ${Q_\emph{\emph{p}}}$  is a straight line $D$  parallel to axis  $\textbf{x}$. Thus, the line $D$   is an image of the inflection circle of the  generalized complex plane by means of an inversion at the rotation center  $I$, (Figure 3,4,5).

	Since the points $Q_\emph{\emph{p}}^1,$   $Q_\emph{\emph{p}}^2$, $Q_\emph{\emph{p}}^3$  are linear, the vectors $\left( {\textbf{IQ}_\emph{\emph{p}}^1 - \textbf{IQ}_\emph{\emph{p}}^2} \right)$  and $\left( {\textbf{IQ}_\emph{\emph{p}}^2 - \textbf{IQ}_\emph{\emph{p}}^3} \right)$  are linearly dependent. So, the  $\emph{\emph{p}}-$cross product of these vectors is
$$\left( {\textbf{IQ}_\emph{\emph{p}}^1 - \textbf{IQ}_\emph{\emph{p}}^2} \right){ \wedge _\emph{\emph{p}}}\left( {\textbf{IQ}_\emph{\emph{p}}^2 - \textbf{IQ}_\emph{\emph{p}}^3} \right) = 0.$$
Thus, we can write
\[\left( {\textbf{IQ}_\emph{\emph{p}}^1{ \wedge _\emph{\emph{p}}}\textbf{IQ}_\emph{\emph{p}}^2} \right) + \left( {\textbf{IQ}_\emph{\emph{p}}^3{ \wedge _\emph{\emph{p}}}\textbf{IQ}_\emph{\emph{p}}^1} \right) + \left( {\textbf{IQ}_\emph{\emph{p}}^2{ \wedge _\emph{\emph{p}}}\textbf{IQ}_\emph{\emph{p}}^3} \right) =0.\]

If we consider the equation (20), we can easily write

\[\frac{1}{{\rho _1^*\rho _2^*}}\left( {\textbf{X}_\emph{\emph{p}}^1{ \wedge _\emph{\emph{p}}}\textbf{X}_\emph{\emph{p}}^2} \right) + \frac{1}{{\rho _1^*\rho _3^*}}\left( {\textbf{X}_\emph{\emph{p}}^3{ \wedge _\emph{\emph{p}}}\textbf{X}_\emph{\emph{p}}^1} \right) + \frac{1}{{\rho _2^*\rho _3^*}}\left( {\textbf{X}_\emph{\emph{p}}^2{ \wedge _\emph{\emph{p}}}\textbf{X}_\emph{\emph{p}}^3} \right) = 0.\]

Considering that $\rho _1^*\rho _2^*\rho _3^* \ne 0$, the last equation and the equation (5), we obtain that

\begin{equation}
\rho _1^*\sin \emph{\emph{p}}\theta _\emph{\emph{p}}^{23} + \rho _2^*\sin \emph{\emph{p}}\theta _\emph{\emph{p}}^{31} + \rho _3^*\sin \emph{\emph{p}}\theta _\emph{\emph{p}}^{12} = 0
\end{equation}
where $\frac{1}{{\rho _k^*}} = \frac{1}{{{\rho _k}}} - \frac{1}{{{{\rho '}_k}}}$  and  $\theta _\emph{\emph{p}}^{lm}$ are the  $\emph{\emph{p}}-$rotation angles between $\textbf{X}_\emph{\emph{p}}^l$  and $\textbf{X}_\emph{\emph{p}}^m$  for the permutations of the indices  $k,l,m = 1,2,3{\kern 1pt} {\kern 1pt} {\kern 1pt} {\kern 1pt} ;{\kern 1pt} {\kern 1pt} {\kern 1pt} 2,3,1{\kern 1pt} {\kern 1pt} {\kern 1pt} ;{\kern 1pt} {\kern 1pt} {\kern 1pt} 3,1,2$.\\

	This formula is called of the generalized Bobillier formula in $C_\emph{\emph{p}}$. The special cases of the generalized Bobillier formula with respect to the sign of real number $\emph{\emph{p}}$  are as follows;

\textbf{Case 1.} If we take  $\emph{\emph{p}} < 0$, the  generalized complex number system  ${C_\emph{\emph{p}}}$ corresponds to the elliptical complex number system. From the equation (6), the generalized Bobillier formula is equal to
\[\rho _1^*\sin \emph{\emph{p}}\left( {\theta _\emph{\emph{p}}^{23}\sqrt {\left| \emph{\emph{p}} \right|} } \right) + \rho _2^*\sin \emph{\emph{p}}\left( {\theta _\emph{\emph{p}}^{31}\sqrt {\left| \emph{\emph{p}} \right|} } \right) + \rho _3^*\sin \emph{\emph{p}}\left( {\theta _\emph{\emph{p}}^{12}\sqrt {\left| \emph{\emph{p}} \right|} } \right) = 0.\]
Especially, if we take  $\emph{\emph{p}} =  - 1$, this formula becomes
\[\rho _1^*\sin {\theta ^{23}} + \rho _2^*\sin {\theta ^{31}} + \rho _3^*\sin {\theta ^{12}} = 0.\]
Under the circumstances, the generalized Bobillier formula is reduced to the Bobillier formula on the complex plane or Euclidean plane. Figure 3 represents the illustration of the Bobillier construction in the complex plane, \cite{erbay:combob, faye:bobfor}.

\begin{center}
\includegraphics[scale=0.9]{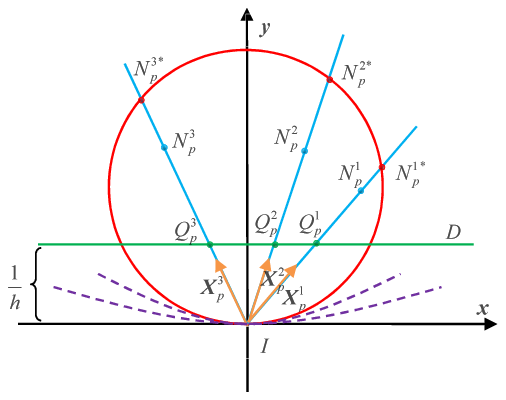}\

\scriptsize
{Figure 3. Inflection Circle in Complex Plane.}
\normalsize
\end{center}

\textbf{Case 2. }If we consider  $\emph{\emph{p}}\, = \,0$, the  generalized complex number system ${C_\emph{\emph{p}}}$  is equal to the parabolic complex number system. By considering the equation (6), we get
\[\rho _1^*\theta^{23} + \rho _2^*\theta ^{31} + \rho _3^*\theta^{12} = 0.\]
This formula is the Bobillier formula for Shear motion in the Galilean plane. The inflection cycle is illustrated as

\begin{center}
\includegraphics[scale=0.9]{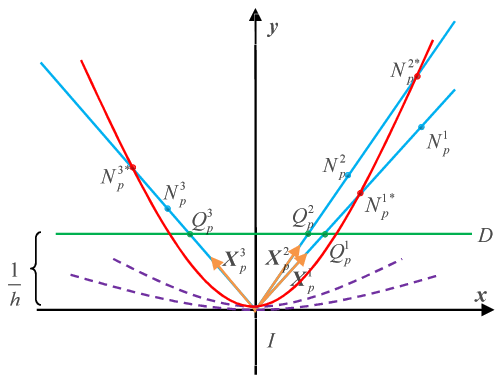}\

\scriptsize
{Figure 4. Inflection Cycle in Galilean Plane.}
\normalsize
\end{center}
for Shear motion in Galilean plane, \cite{guryuc:galbobi}.

\textbf{Case 3. }When  $\emph{\emph{p}} > 0$, the  generalized complex number system ${C_\emph{\emph{p}}}$  is referred to as the hyperbolic complex number system. In addition, considering the equation (6), we obtain that the generalized Bobillier formula is reduced to
\[\rho _1^*\sinh \left( {\theta _\emph{\emph{p}}^{23}\sqrt \emph{\emph{p}} } \right) + \rho _2^*\sinh \left( {\theta _\emph{\emph{p}}^{31}\sqrt \emph{\emph{p}} } \right) + \rho _3^*\sinh \left( {\theta _\emph{\emph{p}}^{12}\sqrt \emph{\emph{p}} } \right) = 0.\]
Especially, if we take  $\emph{\emph{p}} = \,1$, we obtain that the Bobillier formula as
\[\rho _1^*\sinh \theta ^{23} + \rho _2^*\sinh \theta ^{31} + \rho _3^*\sinh \theta ^{12} = 0.\]
Thus, this formula is equal to the Bobillier formula at the Lorentzian plane, \cite{erbay:lobob}. The figure of the (hyperbolic) inflection circle for $\emph{\emph{p}} = \,1$  is as follows

\begin{center}
\includegraphics[scale=0.9]{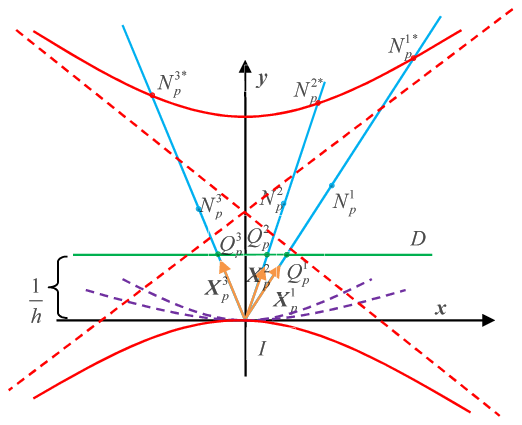}\

\scriptsize
{Figure 5. Inflection Circle in Hyperbolic Plane.}
\normalsize
\end{center}

\section{An Alternative Way for the Generalized Bobillier Formula}

In this section, we will obtain the generalized Bobillier formula in the  generalized complex plane without the use of the generalized Euler-Savary formula. So, the generalized Bobillier formula will be obtained considering the velocities and accelerations of planar  motion in $C_\emph{\emph{p}}$.

Firstly, we calculate the trajectory velocities and accelerations of the points in the moving  generalized complex plane  $K_\emph{\emph{p}}$. Let $\textbf{V}_\emph{\emph{p}}^a\left( {N_\emph{\emph{p}}^1} \right)$  and $\textbf{J}_\emph{\emph{p}}^a\left( {N_\emph{\emph{p}}^1} \right)$  be absolute velocity and acceleration vectors of the point  $N_\emph{\emph{p}}^1$, respectively. Moreover, if  ${w_\emph{\emph{p}}}$ is the angular velocity of  motion  $K_\emph{\emph{p}}/K'_\emph{\emph{p}}$,  ${w_p} = \frac{{d\theta _\emph{\emph{p}}^1}}{{dt}}$ where $\theta _\emph{\emph{p}}^1$  is  $\emph{\emph{p}}-$rotation angle.

	The sliding velocity vector of the inflection point $N{_\emph{\emph{p}}^1}^*$  is perpendicular to the vector which connects the centre to this point. Moreover, the sliding velocity vector is perpendicular to the angular velocity vector. So, we can write

\begin{equation}
\textbf{V}_\emph{\emph{p}}^f\left( {N_\emph{\emph{p}}^1} \right) = {\textbf{w}_\emph{\emph{p}}}{ \wedge _\emph{\emph{p}}}\textbf{IN}_\emph{\emph{p}}^1
\end{equation}
where  ${\textbf{w}_\emph{\emph{p}}} = {w_\emph{\emph{p}}}{\textbf{z}_\emph{\emph{p}}}$ is angular velocity vector of  motion and ${\textbf{z}_\emph{\emph{p}}}$  is the unit vector in the direction of the angular velocity vector.

In addition, from the equation (1) there is the relation between the velocity vectors
\begin{equation}
\textbf{V}_\emph{\emph{p}}^a\left( {N_\emph{\emph{p}}^1} \right) = \textbf{V}_\emph{\emph{p}}^r\left( {{I_\emph{\emph{p}}}} \right) + \textbf{V}_\emph{\emph{p}}^f\left( {N_\emph{\emph{p}}^1} \right)
\end{equation}
where  $\textbf{V}_\emph{\emph{p}}^a$,  $\textbf{V}_\emph{\emph{p}}^f$ and $\textbf{V}_\emph{\emph{p}}^r$  are the absolute, sliding and relative velocity vectors, respectively. If we consider the equations (23) and (24) it is obtained that
\begin{equation}
\textbf{V}_\emph{\emph{p}}^a\left( {N_\emph{\emph{p}}^1} \right) = \textbf{V}_\emph{\emph{p}}^r\left( {{I_\emph{\emph{p}}}} \right) + {\textbf{w}_\emph{\emph{p}}}{ \wedge _\emph{\emph{p}}}\textbf{IN}_\emph{\emph{p}}^1.
\end{equation}
Differentiating the equation (25) with respect to  $t$, we get
\begin{equation}
\textbf{J}_\emph{\emph{p}}^a\left( {N_\emph{\emph{p}}^1} \right) = \textbf{J}_\emph{\emph{p}}^r\left( {{I_\emph{\emph{p}}}} \right) + {\dot w_\emph{\emph{p}}}{\textbf{z}_\emph{\emph{p}}}{ \wedge _\emph{\emph{p}}}\textbf{IN}_\emph{\emph{p}}^1 + w_\emph{\emph{p}}^2\textbf{IN}_\emph{\emph{p}}^1
\end{equation}
where the first term is the path wise tangential acceleration component, the second term is the centripetal component and the third term can be shown to be a pure imaginary component, \cite{saerhurag:eusaveq}. Considering this analysis, the absolute velocity vector of the inflection point is linearly dependent with the absolute acceleration vector of the inflection point since the normal component of acceleration is zero. So we can write
\begin{equation}
\textbf{V}_\emph{\emph{p}}^a\left( {N{{_\emph{\emph{p}}^1}^*}} \right){ \wedge _\emph{\emph{p}}}\textbf{J}_\emph{\emph{p}}^a\left( {N{{_\emph{\emph{p}}^1}^*}} \right) = \textbf{0}.
\end{equation}
Here considering  $\textbf{V}_\emph{\emph{p}}^r = \textbf{0}$  for the inflection point $N{_\emph{\emph{p}}^1}^*$  and using the equations (25) and (26) we can easily find
\[{w_\emph{\emph{p}}}{\left\langle {\textbf{J}_\emph{\emph{p}}^r\left( {{I_\emph{\emph{p}}}} \right),\textbf{IN}_\emph{\emph{p}}^1} \right\rangle _\emph{\emph{p}}} + w_\emph{\emph{p}}^3\rho _1^* = 0\]
and finally it is said that
\begin{equation}
\rho _1^* =  - \frac{{{{\left\langle {\textbf{J}_\emph{\emph{p}}^r\left( {{I_\emph{\emph{p}}}} \right),\textbf{X}_\emph{\emph{p}}^1} \right\rangle }_\emph{\emph{p}}}}}{{w_\emph{\emph{p}}^2}}.
\end{equation}
With similar process for the inflection points $N{_\emph{\emph{p}}^2}^*$  and $N{_\emph{\emph{p}}^3}^*$ the following equations hold;
\begin{equation}
\rho _2^* =  - \frac{{{{\left\langle {\textbf{J}_\emph{\emph{p}}^r\left( {{I_\emph{\emph{p}}}} \right),\textbf{X}_\emph{\emph{p}}^2} \right\rangle }_\emph{\emph{p}}}}}{{w_\emph{\emph{p}}^2}}
\end{equation}
and
\begin{equation}
\rho _3^* =  - \frac{{{{\left\langle {\textbf{J}_\emph{\emph{p}}^r\left( {{I_\emph{\emph{p}}}} \right),\textbf{X}_\emph{\emph{p}}^3} \right\rangle }_\emph{\emph{p}}}}}{{w_\emph{\emph{p}}^2}}.
\end{equation}
On the other hand the linear connection between  $\textbf{X}_\emph{\emph{p}}^1,$ $\textbf{X}_\emph{\emph{p}}^2$  and $\textbf{X}_\emph{\emph{p}}^3$  may be written as follows
\[{\lambda _1}\textbf{X}_\emph{\emph{p}}^1 + {\lambda _2}\textbf{X}_\emph{\emph{p}}^2 + {\lambda _3}\textbf{X}_\emph{\emph{p}}^3 = 0\]
where  ${\lambda _1},{\lambda _3},{\lambda _3} \in R$. From the last equation by successive  $\emph{\emph{p}}-$cross product with $\textbf{X}_\emph{\emph{p}}^1$, $\textbf{X}_\emph{\emph{p}}^2$  and  $\textbf{X}_\emph{\emph{p}}^3$, we get
$${\lambda _1} = \sin \emph{\emph{p}}\theta _\emph{\emph{p}}^{23},{\kern 1pt} {\kern 1pt} {\kern 1pt} {\kern 1pt} {\kern 1pt} {\kern 1pt} {\kern 1pt} {\kern 1pt} {\kern 1pt} {\kern 1pt} {\lambda _2} = \sin \emph{\emph{p}}\theta _\emph{\emph{p}}^{31},{\kern 1pt} {\kern 1pt} {\kern 1pt} {\kern 1pt} {\kern 1pt} {\kern 1pt} {\kern 1pt} {\kern 1pt} {\kern 1pt} {\kern 1pt} {\lambda _3} = \sin \emph{\emph{p}}\theta _\emph{\emph{p}}^{12}$$
where  $\theta _\emph{\emph{p}}^{23} = \theta _\emph{\emph{p}}^3 - \theta _\emph{\emph{p}}^2,$ $\theta _\emph{\emph{p}}^{31} = \theta _\emph{\emph{p}}^1 - \theta _\emph{\emph{p}}^3$   and   $\theta _\emph{\emph{p}}^{12} = \theta _\emph{\emph{p}}^2 - \theta _\emph{\emph{p}}^1.$ So, we obtain that
$$\sin \emph{\emph{p}}\theta _\emph{\emph{p}}^{23}\textbf{X}_\emph{\emph{p}}^1 + \sin \emph{\emph{p}}\theta _\emph{\emph{p}}^{31}\textbf{X}_\emph{\emph{p}}^2 + \sin \emph{\emph{p}}\theta _\emph{\emph{p}}^{12}\textbf{X}_\emph{\emph{p}}^3 = 0.$$

If we make scalar product of the last equation with $\frac{{\textbf{J}_\emph{\emph{p}}^r\left( {{I_\emph{\emph{p}}}} \right)}}{{w_\emph{\emph{p}}^2}}$, the following formula is hold

\begin{equation}
\rho _1^*\sin \emph{\emph{p}}\theta _\emph{\emph{p}}^{23} + \rho _2^*\sin \emph{\emph{p}}\theta _\emph{\emph{p}}^{31} + \rho _3^*\sin \emph{\emph{p}}\theta _\emph{\emph{p}}^{12} = 0.
\end{equation}

This is the formula given in the equation (17). So, this direct way gives us the Bobillier formula without using the generalized Euler-Savary formula.

\section{Conclusion}

The angle between the tangent of pole curve at the instantaneous pole center of coupler with respect to the base of a four bar linkage and of the cranks is equal to the angle between the other crank and the collineation axis. This expression is can be verified by the Bobillier's construction by graphically. This has been a major interest to physicists, mathematicians and engineers. At the same time this problem can be solved by an analytical method called the Bobillier formula which is more practical to use. Various geometric and analytical methods have been developed for the Bobillier formula for Euclidean, Lorentzian and Galilean planar motion. Then the following question can be asked. "Is it possible to give one generalized formula for the Bobillier's construction in generalized complex plane including all planes?" Thus, we find that one formula for all cases of $\emph{\emph{p}}$.  Moreover, we check and interpret this formula for the special cases of $\emph{\emph{p}}$. As a consequence we think that this study would be useful at the disciplines of mathematics, engineering and astronomy.

\end{document}